\documentclass[11pt]{article}
\usepackage[french]{babel}
\usepackage[T1]{fontenc}
\usepackage{graphicx, amsmath, amsthm,amssymb}  
\usepackage{tikz}
\usetikzlibrary{matrix,arrows,calc}

\usepackage{kpfonts}
\usepackage[bookmarks=true,pdfborder={0 0 1},
colorlinks=true,urlcolor=blue,citecolor=purple, 
linkcolor=blue,hypertexnames=false]{hyperref}

\usepackage[margin=10pt,font=small,labelfont=bf,
labelsep=endash]{caption}

\newtheorem{theorem}{Th\'eor\`eme}[section]

\newcommand{\field}[1]{\mathbb{#1}}

\newcommand{\Z}{\field{Z}}    
\newcommand{\R}{\field{R}}    
\newcommand{\E}{\field{E}}    
\newcommand{\fP}{\field{P}}   
\newcommand{\T}{\field{T}}    

\newcommand{\one}{{\mathchoice {\rm 1\mskip-4mu l} {\rm 1\mskip-4mu l}
{\rm 1\mskip-4.5mu l} {\rm 1\mskip-5mu l}}}     

\newcommand{\sB}{\mathscr{B}}
\newcommand{\sI}{\mathscr{I}}
\newcommand{\sL}{\mathscr{L}}

\DeclareMathOperator{\e}{e}          
\DeclareMathOperator{\Hess}{Hess}    
\DeclareMathOperator{\capacity}{cap} 
\DeclareMathOperator{\dd}{d}         
\DeclareMathOperator{\Tr}{Tr}        
\newcommand{\6}[1]{\dd\!#1}          

\newcommand{\eps}{\varepsilon}

\newcommand{\norm}[1]{\left\|#1\right\|} 

\newcommand{\leqs}{\mathbin{\leqslant}}
\newcommand{\geqs}{\mathbin{\geqslant}}

\begin{document}

\title{M\'etastabilit\'e d'EDP stochastiques \\
et d\'eterminants de Fredholm}

\author{Nils Berglund}

\date{15 janvier 2020}

\maketitle

\begin{abstract}
La m\'etastabilit\'e appara\^it lorsqu'un syst\`eme thermodynamique, tel que 
l'eau en surfusion (qui est liquide \`a temp\'erature n\'egative), se 
retrouve du \og mauvais \fg\ c\^ot\'e d'une transition de phase, et reste 
pendant un temps tr\`es long dans un \'etat diff\'erent de son \'etat 
d'\'equilibre. Il existe de nombreux mod\`eles math\'ematiques d\'ecrivant ce 
ph\'enom\`ene, dont des mod\`eles sur r\'eseau \`a dynamique stochastique. Dans 
ce texte, nous allons nous int\'eresser \`a la m\'etastabilit\'e dans des 
\'equations aux d\'eriv\'ees partielles stochastiques (EDPS) paraboliques. 
Certaines de ces \'equations sont mal pos\'ees, et ce n'est que gr\^ace \`a des 
progr\`es tr\`es r\'ecents dans la th\'eorie des EDPS dites singuli\`eres qu'on 
sait construire des solutions, via \`a une proc\'edure de renormalisation. 
L'\'etude de la m\'etastabilit\'e dans ces syst\`emes fait appara\^itre des 
liens inattendus avec la th\'eorie des d\'eterminants spectraux, dont les 
d\'eterminants de Fredholm et de Carleman--Fredholm. 

\smallskip\noindent
\textbf{Article paru dans la \textit{Gazette des Math\'ematiciens}, 
N$^\circ$ 163, Janvier 2020.} 
\end{abstract}

\section{Introduction}

D\'eposez une bouteille d'eau dans le compartiment \`a glace de votre 
r\'efrig\'erateur. Si l'eau est assez pure, en retirant la bouteille apr\`es 
quelques heures, vous trouverez l'eau qu'elle contient encore \`a l'\'etat 
liquide, bien qu'\`a une temp\'erature n\'egative. On dit que l'eau est dans un 
\'etat de surfusion. Agitez la bouteille, et vous verrez l'eau se transformer 
rapidement en glace. 

L'eau en surfusion est un exemple d'\'etat m\'etastable. Dans un tel \'etat, un 
syst\`eme minimise localement un potentiel thermodynamique, tel que son 
\'energie libre, mais pas globalement. La transition vers son \'etat stable 
n\'ecessite de franchir une barri\`ere d'\'energie, ce qui peut prendre beaucoup 
de temps si seules les fluctuations dues \`a l'agitation thermique entrent en 
jeu. Ainsi, la transformation de l'eau en surfusion en glace se fait par 
nucl\'eation, c'est-\`a-dire par l'apparition de cristaux de glace qui croissent 
petit \`a petit.\footnote{Un cristal sph\'erique de rayon $r$ modifie 
l'\'energie du syst\`eme de deux mani\`eres~: d'une part, le fait que la glace 
est plus stable diminue l'\'energie d'une quantit\'e proportionnelle au volume 
du cristal, donc \`a $r^3$; d'autre part, l'interface entre le cristal et l'eau 
environnante augmente l'\'energie d'un terme proportionnel \`a la surface du 
cristal, donc \`a $r^2$. Pour de petites valeurs de $r$, la seconde contribution 
domine la premi\`ere, alors que c'est l'inverse pour $r$ assez grand. Pour cette 
raison, les cristaux de glace croissent tr\`es lentement tant que leur taille 
est plus petite qu'une valeur critique, pour laquelle le terme de volume et le 
terme de surface se compensent.} La pr\'esence d'impuret\'es, ou un apport 
d'\'energie de l'ext\'erieur, peuvent toutefois acc\'el\'erer le processus de 
solidification. 

Il existe de nombreux mod\`eles math\'ematiques d\'ecrivant le ph\'enom\`ene de 
la m\'etastabilit\'e. Les premiers \`a avoir \'et\'e \'etudi\'es sont des 
mod\`eles sur r\'eseau, comme le mod\`ele d'Ising avec une dynamique 
stochastique de type Metropolis--Hastings. On trouvera par exemple 
dans~\cite{denHollander04} un panorama de r\'esultats sur la m\'etastabilit\'e 
dans les syst\`emes dynamiques stochastiques sur r\'eseau. Le m\'etastabilit\'e 
appara\^it toutefois \'egalement dans des syst\`emes continus, tels que les 
\'equations diff\'erentielles stochastiques, que nous allons \'evoquer dans 
la section~\ref{sec:EDS}, ainsi que les EDPs stochastiques que nous aborderons 
dans la section~\ref{sec:AC}. 

\section{Diffusions r\'eversibles}
\label{sec:EDS} 

Le mouvement dans $\R^n$ d'une particule Brownienne de masse $m$, soumise \`a 
une force d\'erivant d'un potentiel $V$, une force de frottement visqueuse, et 
des fluctuations thermiques, peut \^etre d\'ecrit par l'\emph{\'equation de 
Langevin} 
\[
 m\frac{\6^2x_t}{\6t^2} = - \nabla V(x_t) - \gamma\frac{\6x_t}{\6t} 
 + \sigma\,\frac{\6W_t}{\6t}\;,
\] 
o\`u $W_t$ est un mouvement Brownien (voir encart~\ref{sec:MB}), $\gamma$ est un 
coefficient de frottement, et le param\`etre positif $\sigma$ est reli\'e \`a la 
temp\'erature.  Nous supposerons dans la suite que $V:\R^n\to\R$ est un 
potentiel confinant (born\'e inf\'erieurement et tendant vers l'infini assez 
rapidement), et nous sommes int\'eress\'es surtout au cas o\`u $\sigma$ est 
petit. De plus, nous \'ecrirons $\sigma=\sqrt{2\eps}$, afin de simplifier un  
certain nombre d'expressions. 

Lorsque $\eps=0$, si la masse $m$ de la particule est assez petite par rapport 
au coefficient de frottement $\gamma$, la particule s'approche sans osciller 
d'un minimum local de $V$. On dit que son mouvement est \emph{suramorti}. Pour 
$\eps$ quelconque et dans la limite de $m/\gamma$ tr\`es petit, on peut montrer 
qu'apr\`es un changement d'unit\'es, le mouvement de la particule Brownienne est 
d\'ecrit par l'\'equation plus simple du premier ordre
\begin{equation}
\label{eq:EDS} 
 \frac{\6x_t}{\6t} = -\nabla V(x_t) +  \sqrt{2\eps}\,\frac{\6W_t}{\6t}\;, 
\end{equation}
qu'on appelle une \'equation de Langevin suramortie. Math\'ematiquement parlant, 
c'est un exemple d'\emph{\'equation diff\'erentielle stochastique} (EDS), et sa 
solution est aussi appel\'ee une \emph{diffusion}.

Par exemple, en dimension $n=1$, si $V(x) = \frac12 x^2$ 
l'\'equation~\eqref{eq:EDS} devient 
\begin{equation}
\label{eq:OU} 
  \frac{\6x_t}{\6t} = -x_t  + \sqrt{2\eps}\,\frac{\6W_t}{\6t}\;,
\end{equation} 
et d\'ecrit un oscillateur harmonique suramorti soumis \`a un bruit thermique. 
Sa solution est appel\'ee un~\emph{processus d'Ornstein--Uhlenbeck}. 

Une mani\`ere de d\'ecrire les solutions de~\eqref{eq:EDS} est de d\'eterminer 
leurs probabilit\'es de transition $p_t(x,y)$. Celles-ci sont telles que si la 
particule d\'emarre du point $x$ au temps $0$, alors la probabilit\'e $\fP^x 
\{ x_t \in A\}$ de la trouver dans une r\'egion $A$ en un temps $t>0$ s'\'ecrit 
comme 
\[
 \fP^x \left\{ x_t \in A \right\}
 = \int_A p_t(x,y) \6y\;.
\]
On sait que $p_t(x,y)$ satisfait~\emph{l'\'equation de Fokker--Planck} 
\begin{equation}
\label{eq:FPE} 
  \partial_t p_t = 
 \nabla \cdot \left(\nabla V p_t\right) + \eps \Delta p_t 
\end{equation} 
(les op\'erateurs $\Delta$ et $\nabla$ agissant sur la variable $y$).
Le terme $\nabla \cdot (\nabla V p_t)$ a pour effet de transporter $p_t$ d'une 
distance proportionnelle \`a $-\nabla V$, alors que $\eps \Delta p_t$ est un 
terme de diffusion, qui tend \`a \'elargir la distribution de $x_t$. 
Dans le cas du processus d'Ornstein--Uhlenbeck~\eqref{eq:OU}, on peut 
v\'erifier que 
\begin{equation}
\label{eq:OUtrans} 
 p_t(x,y) = \frac{1}{\sqrt{2\pi\eps(1-\e^{-2t})}}
 \exp \biggl\{ - \frac{(y-x\e^{-t})^2}{2\eps(1-\e^{-2t})}\biggr\}\;,
\end{equation} 
c'est-\`a-dire que $x_t$ suit une loi normale d'esp\'erance $x\e^{-t}$ et de 
variance $\eps(1-\e^{-2t})$. Remarquons que lorsque $t$ tend vers l'infini, 
cette loi tend vers une loi normale centr\'ee de variance $\eps$~: plus la 
temp\'erature est faible, plus la variance est petite, et les fluctuations de 
$x_t$ sont moins importantes. 

\begin{figure}[tb]
\begin{center}
\includegraphics[width=7cm]{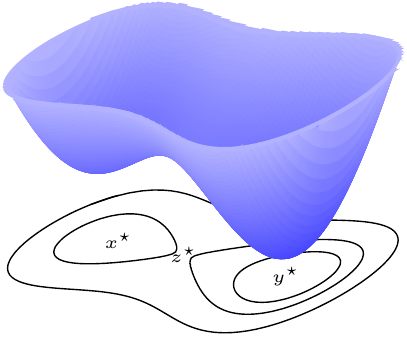}
\end{center}
\vspace{-5mm}
 \caption[]{Un potentiel \`a deux puits. Les minima locaux $x^\star$ et 
$y^\star$ sont s\'epar\'es par un point selle $z^\star$.}
 \label{fig:potentiel}
\end{figure}

Pour des potentiels $V$ g\'en\'eraux, on ne sait pas r\'esoudre l'\'equation de 
Fokker--Planck \eqref{eq:FPE}. Toutefois, on sait que la limite lorsque 
$t\to\infty$ de $p_t(x,y)$ est toujours \'egale \`a 
\[
 \pi(y) = \frac{1}{Z} \e^{-V(y)/\eps} 
\]
o\`u $Z$ est une constante de normalisation telle que l'int\'egrale de $\pi(y)$ 
vaille $1$. En fait,\footnote{L'invariance de $\pi$ suit du fait  que $\pi$ 
appartient au noyau de l'op\'erateur de Fokker--Planck apparaissant au membre de 
droite de l'\'equation \eqref{eq:FPE}, ce qui \'equivaut \`a la condition 
$\eps\nabla\cdot (\e^{-V/\eps}\nabla(\e^{V/\eps}\pi)) = 0$. La relation 
d'\'equilibre d\'etaill\'e~\eqref{eq:eq_detaille} vient du fait que cet 
op\'erateur est auto-adjoint dans l'espace $L^2$ muni du poids $\e^{V/\eps}$.} 
$\pi(y)\6y$ est aussi une mesure de probabilit\'e invariante du processus,  
c'est-\`a-dire que 
\[
 \int_{\R^n} \pi(x) p_t(x,y) \6x = \pi(y)
 \qquad 
 \forall y\in\R^n\;, \forall t>0\;.
\]
Mieux, on sait montrer que la diffusion $(x_t)_{t\geqs0}$ 
est~\emph{r\'eversible} par rapport \`a $\pi$~: ses probabilit\'es de 
transition satisfont la~\emph{condition d'\'equilibre d\'etaill\'e} 
\begin{equation}
\label{eq:eq_detaille} 
 \pi(x) p_t(x,y) = \pi(y) p_t(y,x) 
 \qquad \forall x,y\in\R^n\;, \forall t>0\;. 
\end{equation} 
Cette relation se v\'erifie ais\'ement dans le cas des probabilit\'es de 
transition~\eqref{eq:OUtrans} du processus d'Ornstein--Uhlenbeck.  
Physiquement, elle signifie que si l'on renverse le sens du temps, les 
trajectoires gardent la m\^eme probabilit\'e. Autrement dit, si l'on filmait le 
syst\`eme \`a l'\'equilibre et qu'on passait le film \`a l'envers, on serait 
incapable de d\'etecter une diff\'erence. 

\begin{figure}[tb]
\begin{center}
\includegraphics[width=0.8\textwidth]{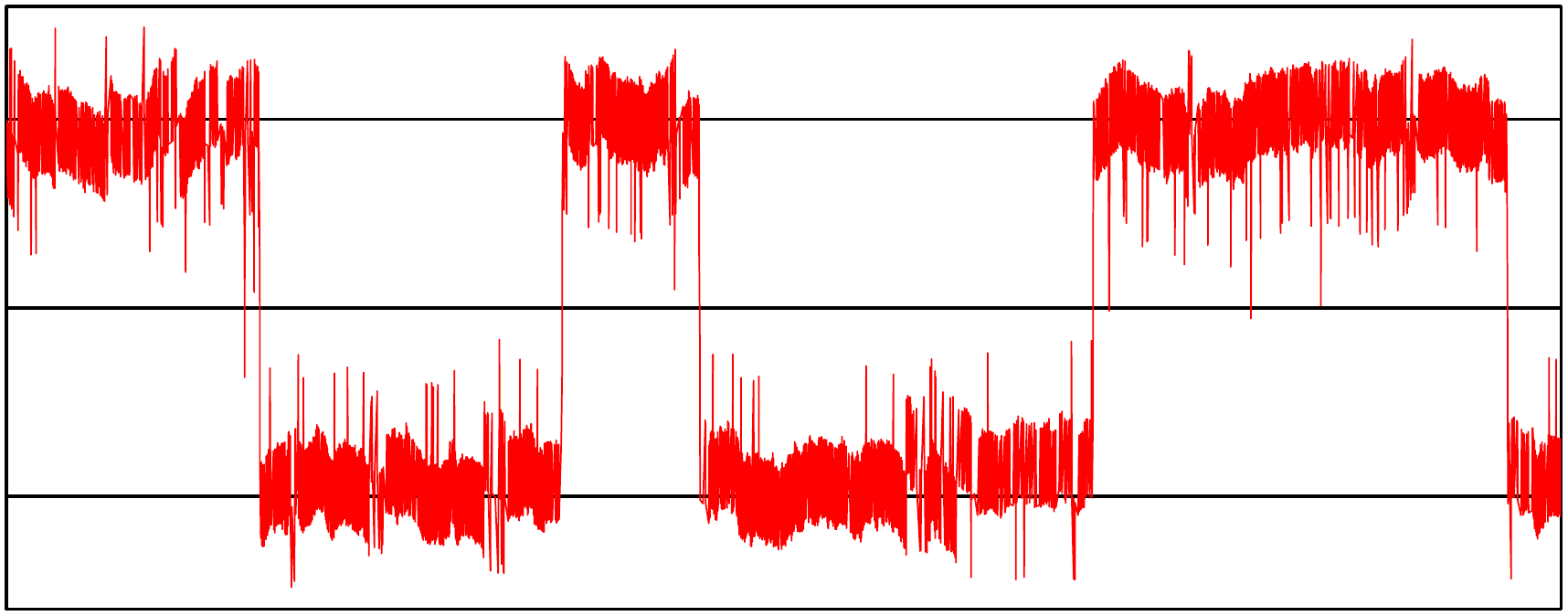}
\end{center}
\vspace{-5mm}
 \caption[]{Une trajectoire $x_t$ de l'EDS~\eqref{eq:EDS} dans un cas de 
dimension $1$, avec le potentiel $V(x) = \frac14x^4 - \frac12x^2$. La 
trajectoire passe la plupart du temps \`a fluctuer autour des deux minima locaux 
$x^\star=-1$ et $y^\star=1$ du potentiel $V$, avec des transitions 
occasionnelles d'un minimum vers l'autre. Dans cette simulation, $\eps$ a 
\'et\'e choisi relativement grand pour que des transitions soient observables 
durant le temps de la simulation.}
 \label{fig:double_puits}
\end{figure}

La m\'etastabilit\'e se manifeste dans le syst\`eme~\eqref{eq:EDS} d\`es que $V$ 
admet plus d'un minimum local. Consid\'erons le cas le plus simple o\`u $V$ est 
un potentiel \`a deux puits, c'est-\`a-dire que $V$ admet exactement deux minima 
locaux $x^\star$ et $y^\star$, ainsi qu'un point selle $z^\star$ 
(Figure~\ref{fig:potentiel}). Les deux minima locaux repr\'esentent deux \'etats 
m\'etastables du syst\`eme, car les solutions de l'EDS~\eqref{eq:EDS} passent 
beaucoup de temps au voisinage de ces points (Figure~\ref{fig:double_puits}). 

La question centrale est alors la suivante. Supposons que la diffusion d\'emarre 
dans le premier minimum local $x^\star$, et soit $\sB_\delta(y^\star)$ une boule 
de petit rayon $\delta$ centr\'ee au second minimum. Quel est le comportement, 
pour $\eps$ petit, du premier temps o\`u $x_t$ visite $\sB_\delta(y^\star)$, 
not\'e $\tau = \inf\{t>0 \mid x_t\in\sB_\delta(y^\star)\}$~? 

\subsection{Loi d'Arrhenius et th\'eorie des grandes d\'eviations}
\label{ssec:LDP} 

Une premi\`ere r\'eponse \`a cette question fut propos\'ee d\`es la fin du XIXe 
si\`ecle par Jacobus van t'Hoff, puis justifi\'ee physiquement par Svante 
Arrhenius~\cite{Arrhenius}~: la valeur moyenne de $\tau$ (son esp\'erance) se 
comporte comme $\e^{[V(z^\star)-V(x^\star)]/\eps}$. Elle est donc 
exponentiellement grande dans la hauteur de la barri\`ere de potentiel entre les 
deux minima locaux de $V$. Lorsque $\eps$ tend vers $0$, le temps de transition 
moyen tend tr\`es rapidement vers l'infini, refl\'etant le fait qu'aucune 
transition n'est possible en l'absence de fluctuation thermique. Inversement, 
lorsque $\eps$ augmente, le temps moyen entre transitions devient de plus en 
plus court. 

Une version rigoureuse de cette loi dite~\emph{d'Arrhenius} peut \^etre 
d\'eduite de la \emph{th\'eorie des grandes d\'eviations}, d\'evelopp\'ee dans 
le contexte des EDS par Mark Freidlin et Alexander Wentzell dans les ann\'ees 
1960--70~\cite{FW}. L'id\'ee de l'approche est la suivante. On fixe un 
intervalle de temps $[0,T]$, et on associe \`a toute trajectoire 
d\'eterministe d\'erivable $\gamma:[0,T]\to\R^n$ la~\emph{fonction taux} 
\begin{equation}
\label{eq:fct_taux} 
 \sI_{[0,T]}(\gamma) = \frac12 \int_0^T \norm{\frac{\6\gamma}{\6t}(t) + \nabla 
V(\gamma(t))}^2 \6t\;. 
\end{equation}
Remarquons que cette fonction est nulle si et seulement si $\gamma(t)$ satisfait 
l'\'equation $\frac{\6\gamma}{\6t}(t)=-\nabla V(\gamma(t))$, 
c'est-\`a-dire~\eqref{eq:EDS} pour $\eps=0$. Sinon, $\sI_{[0,T]}(\gamma)$ est 
strictement positive, et mesure le \og co\^ut \fg\ pour que $x_t$ reste proche 
de $\gamma(t)$. En effet, le principe des grandes d\'eviations pour les 
diffusions affirme que la probabilit\'e que cela arrive est proche (dans un 
sens pr\'ecis) de l'exponentielle de $- \sI_{[0,T]}(\gamma)/(2\eps)$. 

On peut \'egalement estimer la probabilit\'e $p(T) = \fP^{x^\star}\! \left\{ 
\tau \leqs T \right\}$ que la diffusion partant de $x^\star$ atteigne la boule 
$\sB_\delta(y^\star)$ en un temps $T$ au plus. Observons pour cela que pour 
tout $T_1 \in[0,T]$, la fonction taux est sup\'erieure ou \'egale \`a 
$\sI_{[0,T_1]}(\gamma)$, qui peut aussi s'\'ecrire 
\[
 \sI_{[0,T_1]}(\gamma)
 = \frac12 \int_0^{T_1} \norm{\frac{\6\gamma}{\6t}(t) - \nabla V(\gamma(t))}^2 
\6t + 2 \int_0^{T_1} \frac{\6\gamma}{\6t}(t) \cdot \nabla V(\gamma(t)) \6t\;.
\]
Le second terme s'int\`egre et vaut $2[V(\gamma(T_1)) - V(\gamma(0))]$. Comme 
le potentiel le long de toute trajectoire $\gamma$ reliant $x^\star$ \`a 
$\sB_\delta(y^\star)$ atteint au moins la valeur $V(z^\star)$, le principe de 
grandes d\'eviations montre que $p(T)$ est au plus d'ordre 
$\e^{-[V(z^\star)-V(x^\star)]/\eps}$. On peut de plus construire une 
trajectoire de $x^\star$ \`a $\sB_\delta(y^\star)$ de co\^ut $2[V(z^\star) - 
V(x^\star) + R(T)]$ o\`u $R(T)$ est un reste tendant vers $0$ lorsque 
$T\to\infty$.\footnote{Pour $T$ assez grand, on relie des points voisins 
de $x^\star$ et $z^\star$ en un temps $(T-1)/2$ par une trajectoire sur 
laquelle $\smash{\frac{\6\gamma}{\6t}(t)}=+\nabla V(\gamma(t))$, de co\^ut 
proche de $2[V(z^\star) - V(x^\star)]$. Puis on relie un point proche de 
$z^\star$ \`a $\sB_\delta(y^\star)$ en un temps $(T-1)/2$ par une trajectoire 
d\'eterministe de co\^ut nul. Enfin, on utilise le temps $1$ restant pour 
compl\'eter ces deux bouts de trajectoires par des segments de droites, de 
co\^ut n\'egligeable.}
On conclut alors en comparant le processus \`a un processus de Bernoulli, 
effectuant des tentatives ind\'ependantes d'atteindre $\sB_\delta(y^\star)$ 
pendant les intervalles de temps $[kT,(k+1)T]$, chacune avec probabilit\'e de 
succ\`es $p(T)$, dont l'esp\'erance est \'egale \`a $1/p(T)$. Les erreurs faites 
en comparant les deux processus deviennent en effet n\'egligeables dans la 
limite $\eps\to0$.\footnote{L'\'enonc\'e pr\'ecis du r\'esultat est que le 
temps de transition moyen $\E^{x^\star}[\tau]$ satisfait 
$\lim_{\eps\to0}\eps\ln\E^{x^\star}[\tau]=V(z^\star)-V(x^\star)$.}

\subsection{Loi d'Eyring--Kramers et th\'eorie du potentiel}
\label{ssec:EK} 

La \emph{loi d'Eyring--Kramers}, propos\'ee dans les ann\'ees 
1930~\cite{Eyring,Kramers}, est plus pr\'ecise que la loi 
d'Arrhenius,\footnote{La loi d'Eyring--Kramers a en effet \'et\'e propos\'ee 
une trentaine d'ann\'ees avant qu'on ne dispose d'une preuve de la loi 
d'Arrhenius.} 
puisqu'elle d\'ecrit le pr\'efacteur du temps de transition moyen. D\'enotons 
par $\Hess V(x)$ la matrice Hessienne du potentiel $V$ au point $x$, qu'on 
supposera toujours non d\'eg\'en\'er\'ee (c'est-\`a-dire de d\'eterminant non 
nul). Toutes les valeurs propres de la matrice $\Hess V(x^\star)$ sont 
positives, alors que la matrice $\Hess V(z^\star)$ admet une unique valeur 
propre n\'egative, que nous noterons $\lambda_-(z^\star)$.\footnote{En effet, si 
$\Hess V(z^\star)$ admettait plusieurs valeurs propres n\'egatives, on pourrait 
trouver un chemin plus \'economique en termes d'altitude maximale pour aller de 
$x^\star$ \`a $y^\star$. Par exemple, en dimension $2$, les points stationnaires 
de $V$ auxquels la Hessienne admet deux valeurs propres n\'egatives sont des 
maxima locaux de $V$, alors que nous sommes int\'eress\'es aux cols, 
caract\'eris\'es par une valeur propre de chaque signe.}

Dans cette situation, la loi d'Eyring--Kramers affirme que 
\begin{equation}
\label{eq:EK} 
 \E^{x^\star}\! \left[ \tau \right]
 = \frac{2\pi}{|\lambda_-(z^\star)|} 
 \sqrt{\frac{|\det\Hess V(z^\star)|}{\det\Hess V(x^\star)}}
 \e^{[V(z^\star) - V(x^\star)]/\eps}
 \bigl[ 1 + R(\eps) \bigr]\;, 
\end{equation} 
o\`u $R(\eps)$ est un reste tendant vers $0$ dans la limite $\eps\to0$. Il 
existe actuellement plusieurs m\'ethodes permettant de d\'emontrer ce 
r\'esultat. Dans la suite de cette section, nous allons expliquer celle bas\'ee 
sur la th\'eorie du potentiel, d\'evelopp\'ee par Anton Bovier, Michael Eckhoff, 
V\'eronique Gayrard et Markus Klein dans les ann\'ees 2000~\cite{BEGK}, qui se 
pr\^ete \`a une g\'en\'eralisation aux EDPs stochastiques (les lecteurs qui ne 
s'int\'eressent pas \`a ces pr\'ecisions techniques sont invit\'es \`a passer 
directement \`a la section~\ref{sec:AC}).

Fixons deux ensembles disjoints $A, B\subset\R^n$, \`a bord lisse --- pensez \`a 
des voisinages des minima $x^\star$ et $y^\star$ du potentiel $V$. L'observation 
de base est que la \emph{formule de Dynkin} (ou formule d'It\^o pour les temps 
d'arr\^et) permet d'exprimer plusieurs quantit\'es probabilistes int\'eressantes 
comme solutions d'EDPs. Par exemple, la fonction $w_B(x) = \E^x[\tau_B]$, 
donnant l'esp\'erance du temps d'atteinte de $B$ partant de $x$, satisfait 
le~\emph{probl\`eme de Poisson} 
\begin{equation}
\label{eq:Poisson} 
  \begin{cases}
  (\sL w_B)(x) = -1 &\qquad x\in B^c\;, \\
  w_B(x) = 0        &\qquad x\in B\;,
 \end{cases}
\end{equation}
o\`u $\sL$ est l'op\'erateur diff\'erentiel 
\[
 \sL = \eps\Delta - \nabla V \cdot \nabla\;,  
\]
appel\'e \emph{g\'en\'erateur} de la diffusion $(x_t)_{t\geqs0}$ (c'est 
l'adjoint dans $L^2$ de l'op\'erateur de Fokker--Planck apparaissant 
dans~\eqref{eq:FPE}).

La solution de l'\'equation de Poisson~\eqref{eq:Poisson} peut \^etre 
repr\'esent\'ee sous la forme 
\begin{equation}
\label{eq:wB} 
 w_B(x) = -\int_{B^c} G_{B^c}(x,y)\6y\;,
\end{equation} 
o\`u $G_{B^c}$ est la \emph{fonction de Green} associ\'ee \`a $B^c$, 
solution de 
\[
  \begin{cases}
  (\sL G_{B^c})(x,y) = \delta(x-y) &\qquad x\in B^c\;, \\
  G_{B^c}(x,y) = 0        &\qquad x\in B\;.
 \end{cases}
\] 
La r\'eversibilit\'e implique que $G_{B^c}$ satisfait la~\emph{relation 
d'\'equilibre d\'etaill\'e}
\begin{equation}
\label{eq:sym_GB} 
 \e^{-V(x)/\eps} G_{B^c}(x,y)
 = \e^{-V(y)/\eps} G_{B^c}(y,x)
 \qquad \forall x,y\in B^c\;. 
\end{equation} 
Dans le cas $V=0$, la fonction de Green a une interpr\'etation 
\'electrostatique~: $G_{B^c}(x,y)$ est la valeur en $x$ du potentiel 
\'electrique cr\'e\'e par une charge unit\'e plac\'ee en $y$, lorsque le 
domaine $B$ est occup\'e par un conducteur \`a potentiel nul. 

\begin{figure}[tb]
\begin{center}
\begin{tikzpicture}[main node/.style={circle,minimum
size=0.1cm,inner sep=0pt,fill=white,draw}]


\draw[black,semithick] (-1.5,0.5) -- (-1.5,-1) -- (-0.1,-1) -- (-0.1,-1.3) -- 
(-0.1,-0.7);

\draw[black,semithick] (1.5,0.5) -- (1.5,-1) -- (0,-1);

\draw[black,very thick] (0,-0.85) -- (0,-1.15);

\draw[blue!80,thick,fill=blue!35] plot[smooth cycle,tension=0.8]
  coordinates{(-2.3,0) (-1,0) (-1,1.5) (-2.3,2)};

\draw[blue!80,thick,fill=blue!35] plot[smooth cycle,tension=0.8]
  coordinates{(2.3,0) (1,0.5) (1,2) (2.3,2)};
  
\draw[blue!60,thick] plot[smooth, tension=0.8]
  coordinates{(-0.5,-0.5) (-0.3,0.5) (-0.35,1.5) (-0.7,2.5)};
  
\draw[blue!60,thick] plot[smooth, tension=0.8]
  coordinates{(0.7,-0.5) (0.4,0.5) (0.35,1.5) (0.5,2.5)};
  
\node[blue] at (-1.7,1.5) {$A$};
\node[blue] at (1.6,1.6) {$B$};
\node[] at (-1.6,0.5) {$h_{AB}=1$};
\node[] at (1.7,0.6) {$h_{AB}=0$};

\node[blue!60] at (0.8,2.8) {$h_{AB} = $ constante};

\end{tikzpicture}
\vspace{-3mm}
\end{center}
\caption[]{Le potentiel d'\'equilibre $h_{AB}$ d\'ecrit, dans le cas $V=0$, le 
potentiel \'electrique dans un condensateur form\'e de deux conducteurs $A$ et 
$B$, respectivement au potentiel $1$ et $0$.}
\label{fig:hAB}
\end{figure}
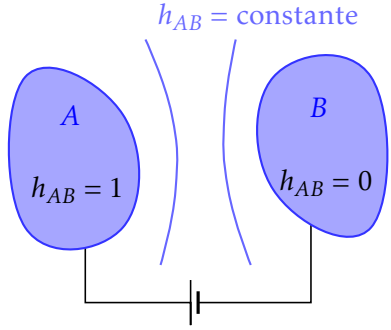

Une seconde quantit\'e importante est le~\emph{potentiel d'\'equilibre} 
$h_{AB}(x) = \fP^x \{\tau_A < \tau_B \}$, aussi appel\'e 
\guillemotleft~committor~\guillemotright~: il donne la probabilit\'e, partant de 
$x$, d'atteindre l'ensemble $A$ avant l'ensemble $B$. C'est une fonction 
${\mathscr L}$-harmonique, qui satisfait le~\emph{probl\`eme de Dirichlet} 
\[
  \begin{cases}
  (\sL h_{AB})(x) = 0  &\qquad x\in (A\cup B)^c\;, \\
  h_{AB}(x) = 1        &\qquad x\in A\;, \\
  h_{AB}(x) = 0        &\qquad x\in B\;.
 \end{cases}
\]
Le potentiel d'\'equilibre admet \'egalement une expression int\'egrale en 
termes de la fonction de Green, \`a savoir  
\begin{equation}
\label{eq:hAB} 
 h_{AB}(x) = -\int_{\partial A} G_{B^c}(x,y) e_{AB}(\6y)\;,
\end{equation} 
o\`u $e_{AB}$ est une mesure concentr\'ee sur  $\partial A$, appel\'ee  
\emph{mesure d'\'equilibre}, d\'efinie par  
\[
 e_{AB}(\6x) = (-{\mathscr L} h_{AB})(\6x)\;.
\]
L'interpr\'etation \'electrostatique de $h_{AB}$ est que c'est le potentiel 
\'electrique dans un con\-densateur, form\'e de deux conducteurs en $A$ et $B$, 
respectivement au potentiel $1$ et $0$ (Figure~\ref{fig:hAB}). 
Enfin, la \emph{capacit\'e} 
\[
 \capacity(A,B) = \int_{\partial A} \e^{-V(x)/\eps} e_{AB}(\6x)
\]
est la constante de normalisation assurant que 
\[
 \nu_{AB}(\6x) = \frac{1}{\capacity(A,B)} \e^{-V(x)/\eps} e_{AB}(\6x)
\] 
soit une mesure de probabilit\'e sur $\partial A$. En \'electrostatique, 
$\capacity(A,B)$ s'interpr\`ete comme la charge totale dans le condensateur 
(qui est bien \'egale \`a la capacit\'e pour une diff\'erence de potentiel 
unit\'e). 

\begin{figure}[tb]
\begin{center}
\begin{tikzpicture}[main node/.style={circle,minimum
size=0.1cm,inner sep=0pt,fill=white,draw},scale=0.9]


\draw[semithick] (0,0) -- (3,2);
\draw[semithick] (0,0) -- (4,0.5);
\draw[] (2.7,1.8) -- (2.9,1.5) -- (3.2,1.7);
\draw[semithick,blue!75] (5,-1) -- (2,3.5);

\node[main node] at (0,0) {};
\node[main node] at (3,2) {};
\node[main node] at (4,0.5) {};

\node[] at (-0.3,-0.2) {$0$};
\node[] at (3.5, 2.2) {$h_{AB}$};
\node[] at (4.4, 0.6) {$f$};
\node[blue!75] at (5.4, -0.8) {${\mathscr H}_{AB}$};

\end{tikzpicture}
\vspace{-3mm}
\end{center}
\caption[]{Le principe de Dirichlet affirme que la capacit\'e minimise la 
distance \`a l'origine, mesur\'ee par la forme de Dirichlet, parmi l'ensemble  
${\mathscr H}_{AB}$ des fonctions $h$ valant $1$ en $A$ et $0$ en $B$. Cela 
d\'ecoule du fait que $\langle f,-{\mathscr L} h_{AB}\rangle_\pi$ est constante 
pour $f\in {\mathscr H}_{AB}$.}
\label{fig:Dirichlet}
\end{figure}
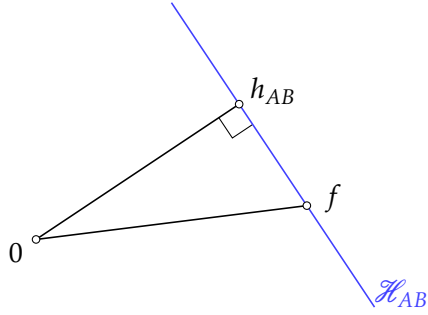

En combinant les expressions~\eqref{eq:wB} de $w_B$ et~\eqref{eq:hAB} de 
$h_{AB}$ avec la relation d'\'equilibre d\'etaill\'e~\eqref{eq:sym_GB} de la 
fonction de Green, on obtient la relation 
\begin{equation}
 \label{eq:magic}
 \int_{\partial A} \E^x[\tau_B]\nu_{AB}(\6x) 
 = \frac{1}{\capacity(A,B)} \int_{B^c} \e^{-V(x)/\eps} h_{AB}(x) \6x\;,
\end{equation} 
qui est essentielle pour l'approche par la th\'eorie du potentiel. En prenant 
pour $A$ une petite boule centr\'ee en $x^\star$, on peut en effet montrer (soit 
\`a l'aide d'in\'egalit\'es de Harnack, soit par un argument de couplage) que 
$\E^x[\tau_B]$ varie tr\`es peu sur $\partial A$. Le membre de gauche 
de~\eqref{eq:magic} est donc proche de l'esp\'erance cherch\'ee 
$\E^{x^\star}[\tau_B]$. Quant au membre de droite, on commence par observer que 
si $B$ est une petite boule centr\'ee en $y^\star$, alors $h_{AB}$ est proche de 
$1$ dans le bassin d'attraction de $x^\star$, et exponentiellement petite (dans 
un sens qu'on contr\^ole) dans le bassin de $y^\star$. La m\'ethode de Laplace 
permet alors de montrer que 
\begin{equation}
\label{eq:EK_integrale} 
 \int_{B^c} \e^{-V(x)/\eps} h_{AB}(x) \6x 
 \simeq \sqrt{\frac{(2\pi\eps)^n}{\det\Hess V(x^\star)}} 
\e^{-V(x^\star)/\eps}\;. 
\end{equation}
Il reste donc \`a estimer la capacit\'e. Ceci peut \^etre fait \`a l'aide de 
principes variationnels. La \emph{forme de Dirichlet} est la forme quadratique  
associ\'ee au g\'en\'erateur, qui peut s'\'ecrire, \`a l'aide d'une 
int\'egration par parties (identit\'e de Green) comme 
\[
  {\mathscr E}(f,f) = \langle f,-{\mathscr L} f\rangle_\pi = \eps \int_{\R^n} 
\e^{-V(x)/\eps} \norm{\nabla f(x)}^2 \6x\;,
\]
o\`u $\langle f,g\rangle_\pi$ est le produit scalaire avec poids $\pi(x)$. Le 
\emph{principe de Dirichlet} affirme que la capacit\'e $\capacity(A,B)$ est 
\'egale \`a la borne inf\'erieure de la forme de Dirichlet sur toutes les 
fonctions valant $1$ en $A$ et $0$ en $B$, et que cet infimum est atteint pour 
$f=h_{AB}$. C'est une cons\'equence directe du fait que $\langle f,-{\mathscr L} 
h_{AB}\rangle_\pi = \capacity(A,B)$ pour tous les $f$ satisfaisant ces m\^emes 
conditions aux bords, et de l'in\'egalit\'e de Cauchy--Schwarz (voir 
Figure~\ref{fig:Dirichlet}). En \'electrostatique, la forme de Dirichlet 
s'interpr\`ete comme l'\'energie \'electrostatique du condensateur, qui est 
effectivement minimale dans l'\'etat d'\'e\-quilibre.  

Une borne inf\'erieure \`a la capacit\'e peut \^etre obtenue \`a l'aide 
du~\emph{principe de Thomson}. Pour un champ de vecteurs $\varphi:\R^n\to\R^n$, 
on d\'efinit la forme quadratique 
\[
 {\mathscr D}(\varphi,\varphi) = \frac{1}{\eps} \int_{(A\cup B)^c} 
\e^{V(x)/\eps} \norm{\varphi(x)}^2\6x\;.
\] 
Le principe de Thomson affirme que l'inverse de la capacit\'e est 
l'infimum de ${\mathscr D}$ sur tous les champs de vecteurs de divergence nulle, 
et dont le flux sur $\partial A$ est \'egal \`a $1$. 

En choisissant des fonctions tests ad\'equates pour les deux principes 
variationnels (qu'on devine en s'inspirant du cas de la dimension $1$, 
qui peut \^etre r\'esolu explicitement), on trouve 
\[
 \capacity(A,B) \simeq  
 \frac{|\lambda_-(z^\star)|} {2\pi}
 \sqrt{\frac{(2\pi\eps)^n}{|\det\Hess V(z^\star)|}}
 \e^{-V(z^\star)/\eps}
\]
En prenant le rapport entre~\eqref{eq:EK_integrale} et cette derni\`ere 
expression, on obtient bien la formule d'Eyring--Kramers~\eqref{eq:EK}.

\section{M\'etastabilit\'e pour l'\'equation d'Allen--Cahn}
\label{sec:AC} 

Notre objectif est maintenant de quantifier la m\'etastabilit\'e, de mani\`ere 
similaire aux EDS, pour des \'equations aux d\'eriv\'ees partielles 
stochastiques (EDPS). Nous allons consid\'erer ici l'\'equation 
d'Allen--Cahn
\begin{equation}
 \label{eq:AC} 
  \partial_t \phi = \Delta\phi + \phi - \phi^3 + \sqrt{2\eps}\xi\;,
\end{equation}
qui est un mod\`ele simple de s\'eparation de phases, dans un m\'elange de glace 
et d'eau par exemple, ou encore dans un alliage. C'est \'egalement l'une 
des EDPS les plus simples pr\'esentant un comportement m\'etastable. 

L'inconnue est un champ $\phi(t,x)$, o\`u la variable spatiale $x$ appartient au 
tore $\T^d_L = (\R/L\Z)^d$ de taille $L$ (on pourrait travailler avec le tore 
unit\'e quitte \`a introduire un param\`etre de viscosit\'e devant le 
Laplacien). Le terme $\xi$ d\'enote ce que l'on appelle un \emph{bruit blanc 
espace-temps}. Intuitivement, $\xi$ repr\'esente un bruit Brownien agissant de 
mani\`ere ind\'ependante en tout point de l'espace, ce qui se traduit par la 
relation
\begin{equation}
 \label{eq:xi_cov} 
 \E[\xi(t,x)\xi(s,y)] = \delta(t-s)\delta(x-y)\;.
\end{equation} 
La d\'efinition math\'ematique de $\xi$ est que c'est une distribution 
al\'eatoire gaussienne, centr\'ee, et de 
covariance donn\'ee par
\[
 \E \! \left[ \langle\xi,\varphi_1 \rangle \langle\xi,\varphi_2 \rangle \right]
 = \langle \varphi_1, \varphi_2 \rangle_{L^2}
\]
pour tout couple de fonctions tests $\varphi_1, \varphi_2\in L^2$. 
En effet, en rempla\c cant formellement les fonctions test par des 
distributions de Dirac, on retrouve la relation~\eqref{eq:xi_cov}. De plus, si 
$\varphi_T(t,x)$ vaut $1$ si $t\in[0,T]$ et $x$ appartient \`a un ensemble 
$A\subset\T_L^d$, et $0$ sinon, alors $W_T=\langle\xi,\varphi_T\rangle$ est un 
mouvement Brownien. 

On peut consid\'erer~\eqref{eq:AC} comme un analogue en dimension infinie de 
la diffusion gradient~\eqref{eq:EDS}, pour le potentiel 
\begin{equation}
\label{eq:pot} 
V(\phi) = \int_{\T_L^d} \biggl( \frac12 \Vert\nabla\phi(x)\Vert^2 - \frac12 
\phi(x)^2 + \frac14 \phi(x)^4 \biggr) \6x\;.
\end{equation} 
En effet, pour toute fonction p\'eriodique $\psi$, la d\'eriv\'ee 
de G\^ateaux de $V$ dans la direction~$\psi$ vaut 
\[
 \lim_{h\to 0} \frac{V(\phi + h\psi) - V(\phi)}{h}
 = \int_{\T_L^d} \bigl( \nabla\phi(x) \cdot \nabla\psi(x) -  
\phi(x)\psi(x) + \phi(x)^3\psi(x) \bigr) \6x\;.
\]
Une int\'egration par parties du terme en $\nabla\phi\cdot\nabla\psi$ montre que 
cette d\'eriv\'ee est pr\'ecis\'ement \'egale au produit scalaire $-\langle 
\Delta\phi + \phi - \phi^3,\psi\rangle_{L^2}$ du membre de droite 
de~\eqref{eq:AC} avec $\psi$, chang\'e de signe. 

Les solutions stationnaires de~\eqref{eq:AC} dans le cas d\'eterministe $\eps=0$ 
sont les points critiques de $V$. Parmi ces points critiques, il n'en existe que 
deux qui sont des minima locaux, et jouent donc le m\^eme r\^ole que $x^\star$ 
et $y^\star$ dans le cas des diffusions~: ce sont les solutions identiquement 
\'egales \`a $\pm1$, que nous noterons $\phi^\star_\pm$. Si $\phi$ repr\'esente 
un m\'elange d'eau et de glace, alors $\phi^\star_-$ et $\phi^\star_+$ 
repr\'esentent respectivement la glace pure et l'eau pure. Selon la  taille $L$ 
du domaine, il existe un ou plusieurs autres points critiques. Dans ce qui suit, 
pour simplifier, nous allons nous concentrer sur le cas $L<2\pi$. Alors, 
l'unique autre point critique est la fonction identiquement nulle, que nous 
noterons $\phi^\star_{\mathrm{trans}}$ car il s'agit de l'\'etat de transition 
pour aller de $\phi^\star_-$ \`a $\phi^\star_+$. Il joue le m\^eme r\^ole que 
$z^\star$ dans le cas des diffusions. 

\begin{figure}[tb]
\begin{center}
\includegraphics[width=7cm]{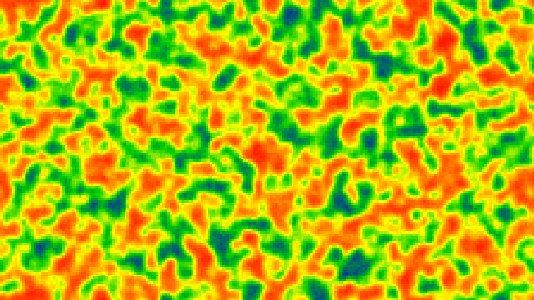}
\hspace{3mm}
\includegraphics[width=7cm]{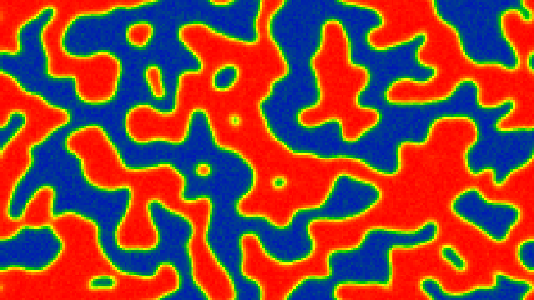} \\
\vspace{2mm}
\includegraphics[width=7cm]{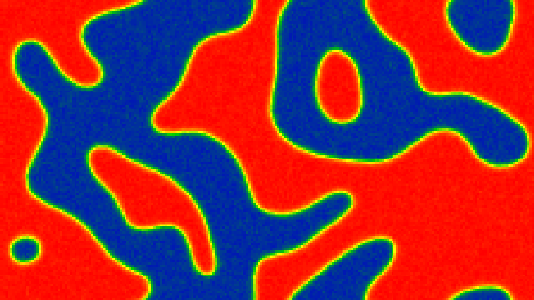}
\hspace{3mm}
\includegraphics[width=7cm]{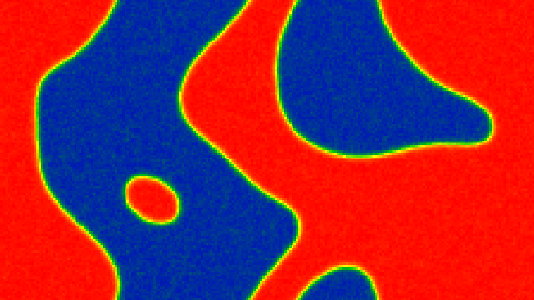}
\end{center}
\vspace{-1mm}
 \caption[]{\'Evolution temporelle d'une solution de l'\'equation 
d'Allen--Cahn stochastique sur un tore de dimension $2$, illustrant le 
ph\'enom\`ene de d\'ecomposition spinodale, ou s\'eparation lente des phases. 
Les couleurs rouge et bleue repr\'esentent, respectivement, des r\'egions o\`u 
le champ $\phi$ est proche de $1$ et de $-1$, alors que le jaune correspond \`a 
$\phi$ proche de $0$. L'effet principal du bruit sur ces simulations est de 
rendre les phases rouge et bleue l\'eg\`erement granulaires. Les interfaces 
entre les deux phases restent relativement lisses, ce qui est d\^u \`a l'effet 
r\'egularisant du Laplacien.}
 \label{fig:Allen-Cahn}
\end{figure}

La Figure~\ref{fig:Allen-Cahn} donne un exemple d'\'evolution temporelle d'une 
solution de~\eqref{eq:AC} en dimension $2$.\footnote{On trouvera des animations 
sur les pages \href{https://www.idpoisson.fr/berglund/simchain.html}{\texttt{ 
https://www.idpoisson.fr/berglund/simchain.html}} pour la dimension $1$, et 
\href{https://www.idpoisson.fr/berglund/simac.html} 
{\texttt{https://www.idpoisson.fr/berglund/simac.html}} pour la dimension $2$. 
Voir aussi la page YouTube \href{http://tinyurl.com/q43b6lf} 
{\texttt{http://tinyurl.com/q43b6lf}}. Par ailleurs, on trouvera des 
simulations interactives aux adresses 
\href{https://experiences.math.cnrs.fr/Equation-aux-Derivees-Partielles.html}
{\texttt{%
https://experiences.math.cnrs.fr/Equation-aux-Derivees-Partielles.html}} et
\href{https://experiences.math.cnrs.fr/EQuation-aux-Derivees-Partielles-69.html}
{\texttt{%
https://experiences.math.cnrs.fr/EQuation-aux-Derivees-Partielles-69.html}
}.} Celle-ci pr\'esente un ph\'enom\`ene de s\'eparation progressive des phases 
(en physique du solide, par exemple pour la s\'eparation des phases d'un 
alliage, on parle aussi de~\emph{d\'ecomposition spinodale}), qui correspond \`a 
une convergence assez lente vers l'un des \'etats stables $\phi^\star_\pm$. Une 
diff\'erence par rapport au cas de la diffusion de dimension $1$ repr\'esent\'e 
dans la Figure~\ref{fig:double_puits} est que l'on observe pendant longtemps un 
m\'elange des deux phases. Ce n'est qu'au bout d'un temps tr\`es long 
(d\'epassant ce qui est montr\'e sur la Figure~\ref{fig:Allen-Cahn}) que le 
syst\`eme s'approche d'une phase pure, soit bleue, soit rouge. Cela est d\^u au 
fait que la condition initiale, al\'eatoire \`a moyenne nulle, incite le champ 
\`a s'approcher d'abord du point selle $\phi^\star_{\mathrm{trans}}$ (qui est 
\'egalement de moyenne nulle), avant d'\^etre finalement attir\'e par 
$\phi^\star_-$ ou $\phi^\star_+$. Comme, de plus, on est en dimension infinie, 
le syst\`eme a beaucoup de \og place \fg\ lui permettant d'\'evoluer avant de 
converger vers un \'equilibre. 

Si, contrairement \`a ce qui est montr\'e dans la Figure~\ref{fig:Allen-Cahn}, 
on d\'emarrait la simulation dans l'une des phases pures, par exemple 
$\phi^\star_-$, on verrait le syst\`eme rester tr\`es longtemps pr\`es de cet 
\'etat, avant de faire une transition vers l'autre phase, c'est-\`a-dire 
$\phi^\star_+$. Puis, apr\`es un autre intervalle de temps tr\`es long, on le 
verrait retourner pr\`es de l'\'etat initial, et ainsi de suite. La valeur 
moyenne du champ se comporterait donc un peu comme dans la 
Figure~\ref{fig:double_puits}. On a donc bien affaire \`a un ph\'enom\`ene de 
m\'etastabilit\'e. Une question naturelle qui se pose ici pour $\eps>0$ est la 
suivante : si l'on d\'emarre avec une condition initiale proche de 
$\phi^\star_-$, quelle est l'asymptotique pr\'ecise du temps n\'ecessaire \`a 
atteindre un petit voisinage (dans une norme appropri\'ee) de la solution 
$\phi^\star_+$ ? 

\subsection{Dimension \texorpdfstring{$1$}{1} : d\'eterminants de Fredholm}

Dans le cas de la dimension $d=1$, William Faris et Giovanni Jona-Lasinio ont 
\'etabli dans~\cite{Faris_JonaLasinio82} un principe de grandes d\'eviations, 
ayant pour fonction taux (comparer \`a l'expression~\eqref{eq:fct_taux} de la 
fonction taux d'une diffusion)
\begin{equation}
\label{eq:taux_AC_1d} 
  {\mathscr I}_{[0,T]}(\gamma) = \frac12 \int_0^T \int_{\T_L} 
\Bigl[\frac{\partial\gamma}{\partial t} (t,x) - 
\frac{\partial^2\gamma}{\partial x^2}(t,x) - \gamma(t,x) + 
\gamma(t,x)^3\Bigr]^2 \6x\6t\;. 
\end{equation}
Soit $\tau$ le temps de premi\`ere atteinte d'une boule 
$B = \{\phi\colon\Vert\phi-\phi^\star_+\Vert_{L^\infty} < \delta\}$, avec 
$\delta>0$ petit et ind\'ependant de $\eps$. 
Par une m\'ethode tout \`a fait analogue \`a celle discut\'ee dans la 
section~\ref{ssec:LDP}, on obtient que $\tau$ satisfait la loi d'Arrhenius
\[
 \E^{\phi^\star_-}\! \left[ \tau \right]
 \simeq \e^{[V(\phi^\star_{\mathrm{trans}})-V(\phi^\star_-)]/\eps}\;.
\]
Qu'en est-il de la loi d'Eyring--Kramers ? Si nous voulons extrapoler 
l'expression~\eqref{eq:EK} obtenue en dimension finie, il nous faut d'abord 
d\'eterminer l'analogue des matrices Hessiennes de $V$ aux points critiques. 
Une int\'egration par parties montre que le d\'eveloppement limit\'e d'ordre 
$2$ du potentiel autour de $\phi^\star_{\mathrm{trans}} = 0$ s'\'ecrit 
\[
 V(\phi) = \frac12 \langle \phi, [-\Delta-1]\phi \rangle_{L^2} 
 + {\mathcal O} (\phi^4)\;,
\]
ce qui nous permet d'identifier $\Hess V(\phi^\star_{\mathrm{trans}})$ avec la 
forme quadratique $-\Delta-1$. Un argument similaire appliqu\'e en 
$\phi^\star_-$ montre que $\Hess V(\phi^\star_-)$ s'identifie \`a 
$-\Delta+2$.\footnote{Les valeurs $-1$ et $2$ sont les d\'eriv\'ees secondes de 
la fonction $\phi\mapsto\frac14\phi^4-\frac12\phi^2$ en $0$ et en $-1$.} Pris 
s\'epar\'ement, ces deux op\'erateurs n'ont pas de d\'eterminant bien d\'efini. 
Toutefois, leur rapport peut s'\'ecrire 
\begin{equation}
 \det\bigl((-\Delta+2)(-\Delta-1)^{-1}\bigr)
  = \det\bigl(\one + 3(-\Delta-1)^{-1}\bigr)\;. 
\label{eq:Fredholm} 
\end{equation}
Il s'agit d'un \emph{d\'eterminant de Fredholm}, un objet qui permet de 
g\'en\'eraliser le polyn\^ome caract\'eristique \`a des op\'erateurs en 
dimension infinie.\footnote{Les racines non nulles du polyn\^ome 
caract\'eristique $c_M(t)=\det(t\one-M)$ d'une matrice $M$ sont les inverses des 
racines de $\bar{c}_M(s) = \det(\one-sM)$. Le d\'eterminant de Fredholm de 
$-sM$ est l'analogue de $\bar{c}_M(s)$ lorsque $M$ est un op\'erateur 
lin\'eaire de dimension infinie.} Pour voir que ce d\'eterminant converge, 
notons que les valeurs propres $\lambda_k$ de $3(-\Delta-1)^{-1}$ d\'ecroissent 
comme $1/k^2$ pour $k$ grand. Le logarithme du d\'eterminant se comporte donc 
comme la somme des $\ln(1+\lambda_k)$, c'est-\`a-dire la somme des $\lambda_k$, 
ou encore la trace de $3(-\Delta-1)^{-1}$. Le crit\`ere de Riemann nous confirme 
que cette somme converge, on dit que  $3(-\Delta-1)^{-1}$ est \emph{de classe 
trace}. En fait, en utilisant deux identit\'es d'Euler sur les produits infinis, 
on peut obtenir la valeur explicite 
\[
 \det\bigl(\one + 3(-\Delta-1)^{-1}\bigr) 
 = - \frac{\sinh^2(L/\sqrt{2})}{\sin^2(L/2)}\;.
\]
Le th\'eor\`eme suivant est un cas particulier d'un r\'esultat montr\'e 
dans~\cite{BG12a} (et aussi d'un r\'esultat de~\cite{Barret15}, obtenu par 
une approche diff\'erente).

\begin{theorem}
\label{thm:EK-1d} 
Pour $L<2\pi$, on a 
\begin{equation}
\label{eq:EK-1d} 
 \E^{\phi^\star_-}\! \left[ \tau \right] = 
 \frac{2\pi}{|\lambda_-(\phi^\star_{\mathrm{trans}})|}
\frac{\e^{[V(\phi^\star_{\mathrm{trans}})-V(\phi^*_-)]/\eps}}{\sqrt{
\bigl|\det\bigl(\one + 
3(-\Delta-1)^{-1}\bigr)\bigr|}} [1+R(\eps,\delta)]\;,
\end{equation} 
o\`u $\lambda_-(\phi^\star_{\mathrm{trans}})=-1$ est la plus petite valeur 
propre de $-\Delta-1$, et $R(\eps,\delta)$ est un terme d'erreur convergeant 
vers $0$ lorsque $\eps\to0$. (La vitesse de cette convergence d\'epend de $L$, 
elle devient plus lente lorsque $L$ s'approche de $2\pi$.) 
\end{theorem}

Donnons une esquisse de la d\'emonstration du Th\'eor\`eme~\ref{thm:EK-1d}. La 
premi\`ere \'etape con\-siste en une~\emph{approximation de Galerkin spectrale}. 
Soit $\{e_k\}_{k\in\Z}$ une base de Fourier de $L^2(\T_L)$, et pour un entier 
positif $N$ (appel\'e \emph{param\`etre de coupure ultra-violette}), soit $P_N$ 
la projection sur l'espace ${\mathscr H}_N$ engendr\'e par $\{e_k\}_{|k|\leqs 
N}$. L'\'equation projet\'ee 
\[
  \partial_t \phi_N = \Delta\phi_N + \phi_N - P_N(\phi_N^3) + 
\sqrt{2\eps}P_N\xi
\]
est alors \'equivalente \`a une EDS de dimension finie de la 
forme~\eqref{eq:EDS}, avec $V$ le potentiel~\eqref{eq:pot} restreint \`a 
${\mathscr H}_N$. On peut donc appliquer l'approche par la th\'eorie du 
potentiel discut\'ee dans la section~\ref{ssec:EK}, en prenant garde \`a bien 
g\'erer la d\'ependance des termes d'erreur dans le param\`etre de coupure $N$, 
puis prendre la limite $N\to\infty$.

Une difficult\'e majeure est donc d'obtenir une estimation similaire 
\`a~\eqref{eq:EK-1d} pour l'approximation de Galerkin, avec un terme d'erreur 
$R(\eps,\delta)$ qui ne d\'epende pas de $N$. Une id\'ee cl\'e de la 
d\'emonstration consiste \`a d\'ecomposer le potentiel $V$ en une partie 
quadratique et une partie d'ordre sup\'erieur. Cela permet d'interpr\'eter la 
capacit\'e et l'int\'egrale du membre de droite de la relation~\eqref{eq:magic} 
comme des esp\'erances, sous une mesure Gaussienne, de certaines variables 
al\'eatoires, que l'on peut ensuite estimer \`a l'aide d'arguments 
probabilistes. On trouvera des d\'etails sur ce calcul dans~\cite[Section 
2.7]{B_Sarajevo}.

\subsection{Dimension \texorpdfstring{$2$}{2} : d\'eterminants de 
Carleman--Fredholm}

Int\'eressons-nous maintenant \`a l'\'equation d'Allen--Cahn~\eqref{eq:AC} sur 
le tore de dimension $d=2$. Il s'av\`ere que contrairement au cas $d=1$, 
l'\'equation n'est plus bien pos\'ee~! C'est une cons\'equence du fait que le 
bruit blanc espace-temps est plus singulier en dimension $2$ qu'en dimension 
$1$. Dans~\cite{daPratoDebussche}, Giuseppe Da Prato et Arnaud Debussche ont 
r\'esolu ce probl\`eme par un proc\'ed\'e de~\emph{renormalisation}, inspir\'e 
de la physique quantique des champs. Au lieu de~\eqref{eq:AC}, ils 
consid\`erent, pour $\delta>0$, l'\'equation r\'egularis\'ee 
\begin{equation}
\label{eq:AC_renorm} 
 \partial_t \phi = \Delta\phi + \phi + 3\eps C_\delta\phi - \phi^3 + 
\sqrt{2\eps}\xi^\delta\;. 
\end{equation} 
Ici $\xi^\delta$ est une r\'egularisation du bruit blanc spatio-temporel, 
d\'efinie comme la convolution $\varrho^\delta*\xi$, o\`u 
\[
\varrho^\delta(t,x) = \frac1{\delta^4} 
\varrho\biggl(\frac t{\delta^2},\frac x{\delta}\biggr)\;,
\]
pour une fonction test $\varrho$ d'int\'egrale $1$. Par cons\'equent, 
$\varrho^\delta$ converge vers la distribution de Dirac lorsque $\delta$ tend 
vers $0$. De plus, $C_\delta$ est une \emph{constante de renormalisation} qui 
diverge comme  $\ln(\delta^{-1})$ lorsque $\delta$ tend vers $0$. Comme 
$\xi^\delta$ est une fonction, et non une distribution, l'\'equation dite 
\emph{renormalis\'ee}~\eqref{eq:AC_renorm} admet des solutions pour tout 
$\delta>0$. Da Prato et Debussche ont alors montr\'e que ces solutions 
convergent vers une limite bien d\'efinie lorsque $\delta$ tend vers $0$. 

\`A premi\`ere vue, on pourrait penser que les \'etats d'\'equilibre stables de 
l'\'equation~\eqref{eq:AC_renorm} se trouvent en $\pm\sqrt{1+3\eps C_\delta}$, 
et tendent donc vers l'infini lorsque $\delta$ tend vers $0$ \`a $\eps$ fix\'e. 
En fait, il n'en est rien --- une premi\`ere indication de cela est que Martin 
Hairer et Hendrik Weber ont d\'emontr\'e dans~\cite{HairerWeber} un principe de 
grandes d\'eviations, avec une fonction taux analogue \`a celle du cas 
unidimensionnel (voir~\eqref{eq:taux_AC_1d}). Le point \`a noter est que comme 
en dimension $1$, cette fonction taux ne fait pas appara\^itre de contre-terme 
de renormalisation. On en d\'eduit la loi d'Arrhenius 
\[
 \E^{\phi^\star_-}\! \left[ \tau \right]
 \simeq \e^{[V(\phi^\star_{\mathrm{trans}})-V(\phi^\star_-)]/\eps}\;, 
\]
o\`u $V$ est le potentiel~\eqref{eq:pot}, ind\'ependant de tout terme de 
renormalisation. Comme avant, $\tau$ est bien le temps de transition entre les 
\'equilibres $\phi^\star_-$ et $\phi^\star_+$, situ\'es en $\pm1$. On 
peut interpr\'eter ce r\'esultat comme indiquant que le \emph{contre-terme} 
$3\eps C_\delta\phi$ sert uniquement \`a rendre la non-lin\'earit\'e $\phi^3$ 
bien d\'efinie. 

Et pour la loi d'Eyring--Kramers~? Ici, il s'av\`ere que le d\'eterminant de 
Fredholm \eqref{eq:Fredholm} ne converge pas. En effet, $3(-\Delta-1)^{-1}$ 
n'est plus de classe trace en dimension $2$, puisque ses valeurs propres sont 
proportionnelles \`a $1/(k_1^2+k_2^2)$ avec $k_1$ et $k_2$ deux entiers non 
nuls. Or la somme de ces valeurs propres diverge comme la s\'erie harmonique~!

La solution \`a ce probl\`eme consiste tout d'abord \`a travailler, comme en 
dimension $1$, avec une approximation de Galerkin spectrale avec coupure 
ultraviolette $N$. Au lieu de r\'egulariser le bruit blanc espace-temps par 
convolution, on peut \`a nouveau consid\'erer sa projection de Galerkin 
spectrale $\xi_N = P_N\xi$, avec un contre-terme 
\[
3\eps C_N = \frac{3\eps}{L^2}\Tr(P_N(-\Delta-1)^{-1})
\]
qui diverge comme $\ln(N)$ (la constante $C_N$ est la variance du \emph{champ 
libre Gaussien} tronqu\'e\footnote{Pour plus d'informations sur le champ libre 
Gaussien, voir l'article de R\'emi Rhodes dans la Gazette n$^\circ$ 157 
(juillet 2018).}). Le potentiel 
renormalis\'e s'\'ecrit alors 
\[
V_N(\phi) = \int_{\T_L^2} \biggl[\frac12\norm{\nabla\phi(x)}^2 + 
\frac14 \phi(x)^4 - \frac12(1+3\eps C_N) \phi(x)^2\biggr] \6x\;.
\]
Le point crucial est alors de noter que 
\[
  V_N(\phi^\star_{\mathrm{trans}}) - V_N(\phi^\star_-) = \frac{L^2}{4} + 
\frac{3}{2}L^2\eps C_N\;.
\]
Le nouveau terme $\frac{3}{2}L^2\eps C_N$ est pr\'ecis\'ement celui qui 
va faire converger le pr\'efacteur. En effet, la 
formule d'Eyring--Kramers fait intervenir le facteur 
\[
 \det\bigl(\one + 3P_N(-\Delta-1)^{-1}\bigr) \e^{-3\Tr(P_N(-\Delta-1)^{-1})}\;, 
\]
qui admet une limite lorsque $N\to\infty$ (cela suit du fait que son logarithme 
se comporte comme la somme des $1/(k_1^2+k_2^2)^2$). Il s'agit en fait d'une 
r\'egularisation connue du d\'eterminant de Fredholm, appel\'ee aussi 
\emph{d\'eterminant de Carleman--Fredholm}, parfois not\'ee $\det_2(\one + 
3(-\Delta-1)^{-1})$. Contrairement au d\'eterminant de Fredholm, ce 
d\'eterminant modifi\'e est bien d\'efini pour les op\'erateurs dont le carr\'e 
est de classe trace, appel\'es \emph{op\'erateurs de Hilbert--Schmidt}, dont 
$3(-\Delta-1)^{-1}$ fait partie. 

Le th\'eor\`eme suivant combine les r\'esultats 
de~\cite{Berglund_DiGesu_Weber_16} et~\cite{Tsatsoulis_Weber_18}. 

\begin{theorem}
\label{thm:EK-2d} 
Soit $\tau$ le temps d'atteinte d'une boule (dans la norme de Sobolev $H^s$ 
pour un $s<0$), centr\'ee en $\phi^\star_+$. Pour $L<2\pi$, on a
\begin{equation}
\label{eq:EK-2d} 
 \E^{\phi^\star_-}\! \left[ \tau \right]
 =  \frac{2\pi}{|\lambda_-(\phi^\star_{\mathrm{trans}})|}
\frac{\e^{[V(\phi^\star_{\mathrm{trans}})-V(\phi^*_-)]/\eps}}{\sqrt{
\bigl|\det_2\bigl(\one + 3(-\Delta-1)^{-1}\bigr)\bigr|}} 
[1+R(\eps,\delta)]\;,
\end{equation} 
o\`u $\lambda_-(\phi^\star_{\mathrm{trans}})=-1$ est la plus petite valeur 
propre de $-\Delta-1$, et $R(\eps,\delta)$ est un terme d'erreur convergeant 
vers $0$ lorsque $\eps\to0$ (\`a une vitesse d\'ependant de $L$.) 
\end{theorem}

Ce r\'esultat confirme que la renormalisation n'a pas pour effet de d\'eplacer 
les \'etats stationnaires, puisque le th\'eor\`eme s'applique bien 
aux \'etats  $\phi^\star_\pm$ situ\'es en $\pm 1$. En revanche, la proc\'edure 
de renormalisation est n\'ecessaire pour obtenir un pr\'efacteur fini pour 
le temps de transition, puisque le rapport de d\'eterminants spectraux et 
le contre-terme $\frac32L^2\eps C_N$ dans le potentiel se compensent 
exactement.

\section{Quelques probl\`emes ouverts}

Une question naturelle est de savoir s'il existe une loi d'Eyring--Kramers pour 
l'\'equation d'Allen--Cahn en dimension $d=3$ (en dimension $4$, on ne s'attend 
pas \`a l'existence de solutions non triviales \`a cette \'equation). Comme 
montr\'e par Martin Hairer dans le tr\`es remarqu\'e 
article~\cite{Hairer2014}\footnote{On pourra consulter l'article de Fran\c cois 
Delarue dans la Gazette n$^\circ$ 143 (janvier 2015) pour plus de d\'etails sur 
la th\'eorie introduite par Martin Hairer, appel\'ee th\'eorie des 
\emph{structures de r\'egularit\'e}.}, qui lui a valu la M\'edaille Fields en 
2014, la forme de l'\'equation renormalis\'ee est alors
\[
 \partial_t \phi = \Delta\phi + \phi + \bigl[3\eps C_\delta^{(1)} - 9\eps^2 
C_\delta^{(2)}\bigr]\phi - \phi^3 + 
\sqrt{2\eps}\xi^\delta\;,
\]
o\`u $C_\delta^{(1)}$ et $C_\delta^{(2)}$ divergent respectivement comme 
$\delta^{-1}$ et $\ln(\delta^{-1})$. Le premier contre-terme provient de la 
m\^eme proc\'edure de renormalisation qu'en dimension $2$ (appel\'ee 
\emph{renormalisation de Wick}), et ne pose pas de nouvelle difficult\'e par 
rapport au cas $d=2$. En revanche, le second contre-terme est propre \`a la 
dimension $3$, et source de nombreuses difficult\'es. En particulier, 
contrairement au cas $d=2$, la mesure invariante de l'\'equation d'Allen--Cahn 
est singuli\`ere par rapport au champ libre Gaussien. 

On peut toutefois remarquer que $(-\Delta-1)^{-1}$ reste Hilbert-Schmidt en 
dimension $3$. Comme le second contre-terme appara\^it avec un facteur 
$\eps^2$, on s'attend \`a ce qu'une formule d'Eyring--Kramers similaire 
\`a~\eqref{eq:EK-2d} reste valable ici. Avec Ajay Chandra, Giacomo Di Ges\`u et 
Hendrik Weber, nous sommes parvenus \`a \'etablir une partie des bornes 
n\'ecessaires \`a \'etablir ce r\'esultat. Toutefois, pour l'heure la borne 
inf\'erieure sur la capacit\'e nous r\'esiste encore.

Bien entendu, il serait souhaitable d'obtenir des formules d'Eyring--Kramers 
pour d'autres EDPS que celle d'Allen--Cahn. Un exemple est l'\'equation de 
Cahn--Hilliard, qui d\'ecrit la s\'eparation de phases dans des situations o\`u 
le volume total de chaque phase est conserv\'e, comme dans le cas de m\'elanges 
d'eau et d'huile. Toutefois, comme la plupart des mod\`eles math\'ematiques de 
syst\`emes m\'etastables, ces EDPS restent bas\'ees sur une dynamique sur 
r\'eseau~: chaque point du r\'eseau est caract\'eris\'e par son \'etat, mais 
reste fix\'e au m\^eme endroit. C'est un bon mod\`ele pour certains alliages ou 
des mat\'eriaux ferromagn\'etiques, qui ont une structure cristalline, avec des 
atomes ou des spins de diff\'erents types attach\'es \`a chaque site. Pour un 
m\'elange d'eau et de glace, toutefois, il n'y a pas de r\'eseau sous-jacent. 
L'un des grands d\'efis de la th\'eorie de la m\'etastabilit\'e est d'analyser 
des mod\`eles tenant compte du fait que les cristaux de glace peuvent se 
d\'eplacer \`a travers l'eau liquide, pour former des cristaux plus grands par 
agglom\'eration.

\appendix

\section{Encart~: le mouvement Brownien}
\label{sec:MB} 

Le mouvement Brownien est un mod\`ele math\'ematique pour le mouvement 
erratique d'une particule immerg\'ee dans un fluide, sous l'effet des 
collisions avec les mol\'ecules du fluide. Il fut observ\'e pour la 
premi\`ere fois par le naturaliste \'ecossais Robert Brown en 1827, lors de 
l'\'etude au microscope de grains de pollen. 

Les premi\`eres descriptions math\'ematiques du mouvement Brownien furent 
propos\'ees par le  math\'ematicien fran\c cais Louis Bachelier en 1901, pour 
des applications en finance, et par Albert Einstein en 1905. Des variantes de 
leurs approches furent d\'evelopp\'ees par  Marian Smoluchowski en 1906 et par 
Paul Langevin en 1908. Les calculs d'Einstein permirent \`a Jean Perrin 
d'estimer exp\'erimentalement le nombre d'Avogadro en 1909, ce qui lui valut le 
prix Nobel de physique en 1926. 

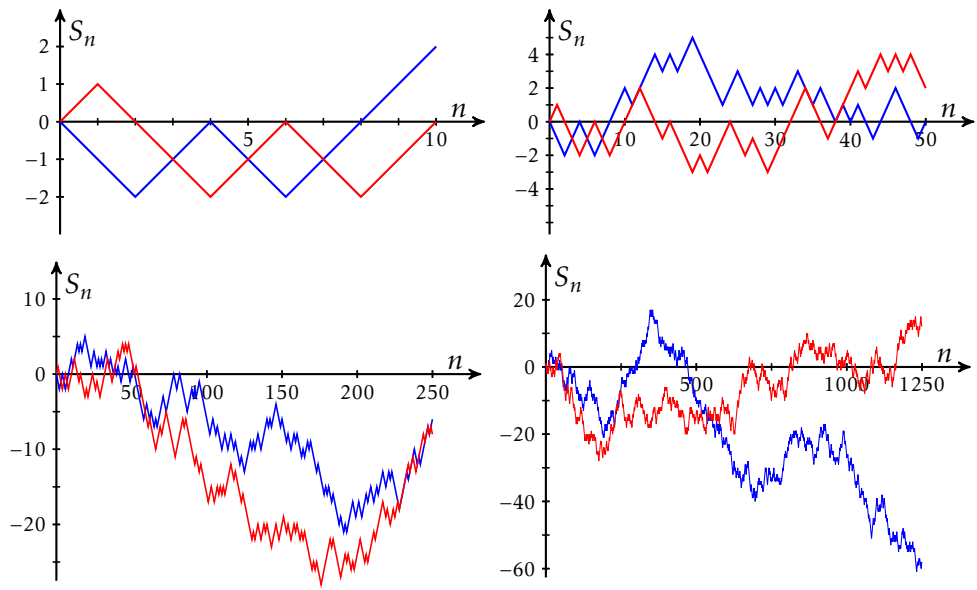
\begin{figure}
\begin{center}
\pgfmathsetmacro{\firstseed}{4}
\pgfmathsetmacro{\secondseed}{12345}
\newcommand*{\firstbmcolor}{blue}
\newcommand*{\secondbmcolor}{red}
\begin{tikzpicture}[>=stealth',main
node/.style={draw,semithick,circle,fill=white,minimum
size=2pt,inner sep=0pt},x=0.5cm,y=0.5cm]



\draw[->,thick] (0,0) -- (11.3,0);
\draw[->,thick] (0,-3) -- (0,3);


\foreach \x in {1,...,10}
\draw[semithick] (\x,-0.1) -- (\x,0.1);
\foreach \x in {5,10}
\draw[semithick] (\x,-0.1) -- node[below] {{\scriptsize $\x$}}
(\x,0.1);

\draw[semithick,white] (-0.1,-2) -- node[left] {{\scriptsize $-20$}}
(0.1,-2);	

\foreach \y in {-2,...,2}
\draw[semithick] (-0.1,\y) -- (0.1,\y);
\foreach \y in {-2,...,2}
\draw[semithick] (-0.1,\y) -- node[left] {{\scriptsize $\y$}}
(0.1,\y);


\pgfmathsetseed{\firstseed}
\draw[thick,\firstbmcolor] (0,0)
\foreach \x in {1,...,10}
{   -- ++(1,{-int(rnd*2)*2+1})
}
;

\pgfmathsetseed{\secondseed}
\draw[thick,\secondbmcolor] (0,0)
\foreach \x in {1,...,10}
{   -- ++(1,{int(rnd*2)*2-1})
}
;


\node[] at (10.6,0.3) {$n$};
\node[] at (0.6,2.35) {$S_n$};
\end{tikzpicture}
\begin{tikzpicture}[>=stealth',main
node/.style={draw,semithick,circle,fill=white,minimum
size=2pt,inner sep=0pt},x=0.5cm,y=0.5cm]


\pgfmathsetmacro{\yscale}{sqrt(0.2)}


\draw[->,thick] (0,0) -- (11.3,0);
\draw[->,thick] (0,-3) -- (0,3);


\foreach \x in {1,...,5}
\draw[semithick] (2*\x,-0.1) -- (2*\x,0.1);
\foreach \x in {1,...,5}
\draw[semithick] (2*\x,-0.1) -- node[below] {{\scriptsize $\x0$}}
(2*\x,0.1);

\draw[semithick,white] (-0.1,-2*\yscale) -- node[left] {{\scriptsize $-20$}}
(0.1,-2*\yscale);	

\foreach \y in {-6,...,5}
\draw[semithick] (-0.1,\y*\yscale) -- (0.1,\y*\yscale);
\foreach \y in {-4,-2,0,2,4}
\draw[semithick] (-0.1,\y*\yscale) -- node[left] {{\scriptsize $\y$}}
(0.1,\y*\yscale);


\pgfmathsetseed{\firstseed}
\draw[thick,\firstbmcolor] (0,0)
\foreach \x in {1,...,50}
{   -- ++(0.2,{(-int(rnd*2)*2+1)*\yscale})
}
;

\pgfmathsetseed{\secondseed}
\draw[thick,\secondbmcolor] (0,0)
\foreach \x in {1,...,50}
{   -- ++(0.2,{(int(rnd*2)*2-1)*\yscale})
}
;


\node[] at (10.6,0.3) {$n$};
\node[] at (0.6,2.35) {$S_n$};
\end{tikzpicture}

\vspace{2mm}
\begin{tikzpicture}[>=stealth',main
node/.style={draw,semithick,circle,fill=white,minimum
size=2pt,inner sep=0pt},x=0.5cm,y=0.5cm]


\pgfmathsetmacro{\yscale}{0.2}


\draw[->,thick] (0,0) -- (11.3,0);
\draw[->,thick] (0,-5.5) -- (0,3);


\foreach \x in {1,...,5}
\draw[semithick] (2*\x,-0.1) -- (2*\x,0.1);

\draw[semithick] (2,-0.1) -- node[below] {{\scriptsize $50$}}
(2,0.1);
\draw[semithick] (4,-0.1) -- node[below] {{\scriptsize $100$}}
(4,0.1);
\draw[semithick] (6,-0.1) -- node[below] {{\scriptsize $150$}}
(6,0.1);
\draw[semithick] (8,-0.1) -- node[below] {{\scriptsize $200$}}
(8,0.1);
\draw[semithick] (10,-0.1) -- node[below] {{\scriptsize $250$}}
(10,0.1);

\foreach \y in {-25,-20,-15,-10,-5,0,5,10}
\draw[semithick] (-0.1,\y*\yscale) -- (0.1,\y*\yscale);
\foreach \y in {-20,-10,0,10}
\draw[semithick] (-0.1,\y*\yscale) -- node[left] {{\scriptsize $\y$}}
(0.1,\y*\yscale);


\pgfmathsetseed{\firstseed}
\draw[semithick,\firstbmcolor] (0,0)
\foreach \x in {1,...,250}
{   -- ++(0.04,{(-int(rnd*2)*2+1)*\yscale})
}
;

\pgfmathsetseed{\secondseed}
\draw[semithick,\secondbmcolor] (0,0)
\foreach \x in {1,...,250}
{   -- ++(0.04,{(int(rnd*2)*2-1)*\yscale})
}
;


\node[] at (10.6,0.3) {$n$};
\node[] at (0.6,2.35) {$S_n$};
\end{tikzpicture}
%
\begin{tikzpicture}[>=stealth',main
node/.style={draw,semithick,circle,fill=white,minimum
size=2pt,inner sep=0pt},x=0.5cm,y=0.5cm]


\pgfmathsetmacro{\yscale}{0.2*sqrt(0.2)}



\draw[->,thick] (0,0) -- (11.3,0);
\draw[->,thick] (0,-5.6) -- (0,3);


\foreach \x in {1,...,5}
\draw[semithick] (2*\x,-0.1) -- (2*\x,0.1);

\draw[semithick] (4,-0.1) -- node[below] {{\scriptsize $500$}}
(4,0.1);
\draw[semithick] (8,-0.1) -- node[below] {{\scriptsize $1000$}}
(8,0.1);
\draw[semithick] (10,-0.1) -- node[below] {{\scriptsize $1250$}}
(10,0.1);

\foreach \y in {-60,-50,-40,-30,-20,-10,0,10,20}
\draw[semithick] (-0.1,\y*\yscale) -- (0.1,\y*\yscale);
\foreach \y in {-60,-40,-20,0,20}
\draw[semithick] (-0.1,\y*\yscale) -- node[left] {{\scriptsize $\y$}}
(0.1,\y*\yscale);


\pgfmathsetseed{\firstseed}
\draw[thin,\firstbmcolor] (0,0)
\foreach \x in {1,...,1250}
{   -- ++(0.008,{(-int(rnd*2)*2+1)*\yscale})
}
;

\pgfmathsetseed{\secondseed}
\draw[thin,\secondbmcolor] (0,0)
\foreach \x in {1,...,1250}
{   -- ++(0.008,{(int(rnd*2)*2-1)*\yscale})
}
;


\node[] at (10.6,0.3) {$n$};
\node[] at (0.6,2.35) {$S_n$};
\end{tikzpicture}
\vspace{-2mm}
\end{center}
\caption[]{Deux r\'ealisations (l'une en rouge, l'autre en bleu) d'une marche 
al\'eatoire sym\'etrique sur $\Z$, vues \`a diff\'erentes \'echelles. D'une 
image \`a la suivante, l'\'echelle horizontale est comprim\'ee d'un facteur 
$5$, alors que l'\'echelle verticale est comprim\'ee d'un facteur $\sqrt{5}$.
}
\label{fig:BM}
\end{figure}

Consid\'erons le cas de la dimension $1$, et supposons que la particule subit 
des collisions r\'eguli\`eres, \`a intervalles de temps $\Delta t$. 
Entre deux collisions successives, la particule se d\'eplace d'une distance 
$\Delta x$, avec probabilit\'e $\frac12$ soit vers la gauche, soit vers la 
droite. Sa position au temps $n\Delta t$ est donc donn\'ee par $S_n\Delta x$, 
o\`u $S_n$ est une suite d'entiers telle que l'incr\'ement $S_{n+1} - S_n$ 
vaille $1$ ou $-1$, chaque fois avec probabilit\'e $\frac12$. On suppose de plus 
que chaque incr\'ement est ind\'ependant de tous les incr\'ements 
pr\'ec\'edents. La suite des $S_n$ s'appelle une marche al\'eatoire 
sym\'etrique sur $\Z$ (voir Figure~\ref{fig:BM}). 

Comme en pratique, les intervalles d'espace et de temps $\Delta x$ et $\Delta 
t$ sont tr\`es petits, il semble pertinent de les faire tendre vers z\'ero, 
afin d'obtenir un objet universel. Il s'av\`ere que cette limite est 
int\'eressante seulement si $\Delta t$ est proportionnel \`a $\Delta x^2$ 
(c'est une cons\'equence du th\'eor\`eme central limite). Cela revient \`a 
poser 
\[
 W_t = \lim_{n\to\infty} \frac{1}{\sqrt{n}} S_{\lfloor nt \rfloor}\;.
\]
Cette d\'efinition s'av\`ere \^etre \'equivalente \`a imposer que pour tout 
$t>s\geqs0$, l'incr\'ement $W_t-W_s$ suive une loi normale, centr\'ee, de 
variance $t-s$, et soit ind\'ependant des valeurs du processus jusqu'au temps 
$s$. 

Norbert Wiener a montr\'e en 1923 que les trajectoires $t\mapsto W_t$ sont 
continues ($W_t$ est d'ailleurs aujourd'hui aussi connu sous le nom de 
\emph{processus de Wiener}). D'autres propri\'et\'es de $W_t$ furent \'etablies 
notamment par Raymond Paley, Antoni Zygmund, et Paul L\'evy. En particulier, on 
sait que les trajectoires du mouvement Brownien ne sont pas diff\'erentiables. 
Ceci pose un probl\`eme pour la d\'efinition de l'EDS~\eqref{eq:EDS}, qu'on 
r\'esout en d\'efinissant ses solutions comme celles de l'\'equation 
int\'egrale
\[
 x_t = x_0 - \int_0^t \nabla V(x_s)\6s + \sqrt{2\eps} W_t\;,
\]
que l'on peut \'etudier par un argument de point fixe. La th\'eorie fut 
g\'en\'eralis\'ee par Kiyoshi It\^o dans les ann\'ees 1940. Son calcul 
stochastique permet de r\'esoudre des variantes de~\eqref{eq:EDS} dans 
lesquelles le terme de bruit est multipli\'e par une fonction de $x$. Certaines 
id\'ees \`a la base du calcul stochastique avaient \'et\'e d\'ecouvertes 
ind\'ependamment par Wolfgang D\"oblin, et envoy\'ees \`a l'Acad\'emie des 
sciences dans un pli cachet\'e qui ne fut ouvert qu'en 2000. 

\bibliographystyle{abbrv}
\bibliography{Meta}

\end{document}